\documentclass[11pt]{article}
\usepackage{latexsym,amsmath,amsfonts}

\newtheorem{theorem}{Theorem}

\newtheorem{lemma}{Lemma}

\newtheorem{remark}{Remark}

\newcommand{\thmref}[1]{Theorem~\ref{thm:#1}} 
\newcommand{\lemref}[1]{Lemma~\ref{lem:#1}} 
\newcommand{\remref}[1]{Remark~\ref{rem:#1}} 
\newcommand{\secref}[1]{Section~\ref{sec:#1}} 
\newcommand{\subref}[1]{Section~\ref{sub:#1}} 
\newcommand{\eqnref}[1]{(\ref{eq:#1})} 

\newcommand\ignore[1]{}

\newcommand{\funcdef}[5]{\begin{array}{rccccc} #1 & : & #2 & \to     & #3 \\
                                                &   & #4 & \mapsto & #5\end{array}}

\def\R{{\bf R}} 
\def\N{{\bf N}} 

\def\periodeq{\mbox{.}}
\def\commaeq{\mbox{,\,}}
\def\semicoloneq{\mbox{;\,}}
\def\suchthat{\,\mbox{:}\,}

\newcommand{\Ex}[1]{\mbox{\bf Ex}\left[#1\right]} 
\newcommand{\Prp}[2]{\mbox{\bf Pr}_{#1}\left(#2\right)} 
\newcommand{\Ind}[1]{\chi_{#1}} 
\newcommand{\Var}[1]{\mbox{\bf Var}\left(#1\right)} 
\def\eqdist{=^d}   
\renewcommand{\Pr}[1]{\mbox{\bf Pr}\left(#1\right)} 

\newcommand{\bigoh}[1]{O\left(#1\right)}
\newcommand{\liloh}[1]{o\left(#1\right)}
\newcommand{\ohmega}[1]{\Omega\left(#1\right)}

\newcommand\QED{\ifhmode\allowbreak\else\nobreak\fi
\quad\nobreak$\Box$\medbreak}
\newcommand{\proofstart}{\par\noindent\sl Proof:\rm\enspace}
\newcommand{\proofend}{\QED\par}
\newenvironment{proof}{\proofstart}{\proofend}

\def\eps{\epsilon}

\def\elead{\mbox{\sf elead}}
\def\periodeq{\mbox{.}}
\def\commaeq{\mbox{,}}

\def\elead{\mbox{\sf ELead}}

\begin{document}

\title{Avoiding defeat in a balls-in-bins process with feedback}
\author{Roberto Oliveira\thanks{IBM T.J. Watson Research Center, Yorktown Heights, NY 10598. \texttt{riolivei@us.ibm.com}, \texttt{rob.oliv@gmail.com}. Done while the author was a Ph.D. student at New York University under the supervision of Joel Spencer. Work funded a CNPq doctoral scholarship.} and Joel Spencer\thanks{Courant Institute of Mathematical Sciences, New York University, New York 10012, NY. Email: \texttt{spencer@cims.nyu.edu}.}} \maketitle
\begin{abstract}Imagine that there are two bins to which balls are added sequentially, and each incoming ball joins a bin with probability proportional to the $p$th power of the number of balls already there. A general result says that if $p>1/2$, there almost surely is some bin that will have more balls than the other at all large enough times, a property that we call {\em eventual leadership}.\\
In this paper, we compute the asymptotics of the probability that
bin $1$ eventually leads when the total initial number of balls $t$
is large and bin $1$ has a fraction $\alpha<1/2$ of the balls; in
fact, this probability is $\exp(c_p(\alpha)t + \bigoh{t^{2/3}})$ for
some smooth, strictly negative function $c_p$. Moreover, we show
that conditioned on this unlikely event, the fraction of balls in
the first bin can be well-approximated by the solution to a certain
ordinary differential equation.
\end{abstract}

\section{Introduction}\label{sec:intro}

Consider a discrete-time process in which there are two bins, to
which balls are added one at a time. Each incoming ball chooses
probabilistically which bin to go to according to the following
rule: if bin $1$ currently has $n_1$ balls and bin $2$ has $n_2$
balls, then the probability that bin $1$ is chosen is
$$\frac{f(n_1)}{f(n_1)+f(n_2)},$$
where $f$ is a fixed positive function. These so-called {\em
balls-in-bins processes with feedback function} $f$\footnote{The
first author's thesis \cite{tese} contains a longer background
discussion of this and related processes.}, which can be generalized
to more than two bins (cf. \secref{prelim} below) were introduced to
the Discrete Mathematics community by Drinea, Frieze and
Mitzenmacher \cite{Drinea02}.  This family of processes was intended
as a model for competition that is mathematically similar to some
so-called preferential attachment models for large networks
\cite{Barabasi99,AlbertSurvey,DrineaEM}.

The authors of \cite{Drinea02} were especially interested in the
case $f(x)=x^p$ with $p>0$ a parameter. In this case, there is a
tendency that {\em the rich get richer}: since $f$ is increasing,
the more balls a bin has, the more likely it is to receive the next
ball. One of the main questions addressed in \cite{Drinea02} is
whether this phenomenon results in effective preponderance by one of
the bins in the long run. They proved that the answer is ``yes" if
$p>1$ and ``no" if $p<1$. That is, if $p>1$ then one of the two bins
will obtain a $1-\liloh{1}$ fraction of all balls in the large-time
limit, whereas if $p<1$ the fractions of balls in the two bins both
tend to $1/2$. The case $p=1$ is the well-known P\'{o}lya Urn, in
which case the limiting number of balls in bin $1$ has a
non-degenerate distribution depending on the initial conditions, so
the result in \cite{Drinea02} seemingly completes the description of
the family of processes given by the choices of $p$.

However, stronger results are available. A paper by Khanin and
Khanin \cite{Khanin01} introduced what amounts to the same process
as a model for neuron growth, and proved that if $p>1$, there almost
surely is some bin that gets {\em all but finitely many balls}, an
event that we call {\em monopoly}. They also show that for
$1/2<p\leq 1$, monopoly has probability $0$, but there almost surely
will be some bin which will lead the process from some finite time
on (we call this {\em eventual leadership}), whereas this cannot
happen if $0<p\leq 1/2$. In fact, the result of \cite{Khanin01}
generalizes to any $f$ with $\min_{x\in\N}f(x)>0$, as shown e.g. in
\cite{Spencer??,tese,EuBrown}.
\begin{theorem}[From \cite{Khanin01,Spencer??,tese,EuBrown}]\label{thm:Khanin}If
$\{I_m\}_{m=0}^{+\infty}$ is a balls-in-bins process and feedback
function $f=f(x)\geq c$ for some $c>0$, then there are three
mutually exclusive possibilities, one of which happens almost surely
irrespective of the initial conditions:
\begin{enumerate}
\item if $\sum_{n\geq 1}f(n)^{-1} < +\infty$, of the
bins receives all but finitely many balls (this is the {\em
monopolistic regime});
\item if $\sum_{n\geq 1}f(n)^{-1} = +\infty$ but
$\sum_{n\geq 1}f(n)^{-2} < +\infty$, monopoly does not happen but
one of the bins has more balls than the other at all large enough
times (this is the {\em eventual leadership regime});
\item if $\sum_{n\geq 1}f(n)^{-2} =+\infty$, the balls alternate in leadership infinitely many times (this is
the {\em almost-balanced regime}).
\end{enumerate}\end{theorem} Notice that the three cases of the Theorem
applied to the $f(x)=x^p$ family correspond to $p>1$, $1/2<p\leq 1$
and $0\leq p<1/2$; in other words, this family of $f$ has {\em phase
transitions} at $p=1$ and $p=1/2$.

The present paper is part of a series of works by the two authors
and by Michael Mitzenmacher in which several more quantitative
aspects of the three regimes are elucidated. We are especially
concerned with the eventual leadership and monopoly regimes, where
there are initially $t$ balls in the system and bin $1$ has a
fraction $\alpha\in(0,1/2)$ of those balls. It is easy to show that
bin $1$ has a positive probability of eventually leading the
process, but this probability should get smaller and smaller as $t$
increases. Thus we ask ourselves two simple questions:

\begin{enumerate}\item {\em How fast does the probability that bin $1$ will
escape its unfavorable initial conditions and eventually lead
converge to $0$ as $t\to +\infty$}? \item {\em What is the typical
behavior of the process, given that bin $1$ does
escape}?\end{enumerate}  Our two main results apply to the case
$f(x)=x^p$, $p>1/2$. We show that the answer to the first question
is ``exponentially small" and compute the exact rate of decay.
Below, let $[t,\alpha]$ denote the pair $(\lceil \alpha
t\rceil,t-\lceil \alpha t\rceil)$.
\begin{theorem}\label{thm:escape_main}Assume that we have a balls-in-bins process with feedback function $f(x)=x^p$, $p>1/2$ (so that the strong eventual leadership condition holds), and let $\elead$ be the event that the first bin eventually leads the process. Then, for all fixed $\alpha\in(0,1/2)$, the limit
\begin{equation}c_p(\alpha) = \lim_{t\to +\infty}\frac{\ln\Prp{[t,\alpha]}{\elead}}{t}\in\R^-\end{equation}
exists, and is a smooth function of $\alpha$ satisfying
$c_p'(\alpha)>0$ on $(0,1/2)$. Moreover, for any $\delta\in(0,1/2)$
there exist $C_\delta\in\R^+$ and $T_\delta\in\N$ such that
\begin{equation}\forall
\alpha\in(\delta,1/2-\delta)\commaeq\forall t\geq
T_\delta\commaeq\;\; e^{c_p(\alpha)\, t - C_\delta\,t^{2/3}} \leq
\Prp{[t,\alpha]}{\elead} \leq e^{c_p(\alpha)\, t +
C_\delta}\periodeq\end{equation}\end{theorem}

The form of \thmref{escape_main} should be compared with that of
Cr\'amer's Theorem \cite{StroockBook}, which estimates the
exponential rate of decay of the probability of large deviations
from the mean of sums of i.i.d. random variables. This analogy also
applies to the proof of \thmref{escape_main} contains computations
with Laplace transforms that resemble those used to prove Cr\'amer's
Theorem. In our case, however, the random variables we consider,
although {\em not} i.i.d., are of a very specific kind.
\thmref{escape_main} is proven in \secref{escape} below.

Question $2.$ tunas out to have a more surprising answer than $1.$.
We will prove that {\em conditioning on bin $1$ escaping almost
determines the behavior of the process}, at least up to the time
when bin $1$ has half of the balls. To state this result precisely,
define
\begin{equation}\funcdef{g_p}{(0,1/2)}{\R}{\alpha}{-\alpha + \frac{\alpha^pe^{c_p'(\alpha)}}{\alpha^pe^{c_p'(\alpha)} + (1-\alpha)^p}}\periodeq\end{equation}
We will show below (cf. \remref{ODE_increasing}) that
$g_p(\alpha)>0$ for all $\alpha\in(0,1/2)$. This implies that the
function $A=A_{\alpha,p}(\cdot)$ solving the following ODE is
increasing.
\begin{equation}\label{eq:ODE_ode}\left\{\begin{array}{llll}\frac{dA}{ds}(s) &=& g_p(A(s))\commaeq & s>0\\
A(0)&=&\alpha & \end{array}\right.\end{equation}

Such a solution is only guaranteed to exist for
$s\in[0,T_{p,\alpha})$, where $T_{p,\alpha}\in\R^+\cup\{+\infty\}$.
Given that $A_{p,\alpha}$ is increasing, $T_{p,\alpha}$ is finite if
and only if $\lim_{s\to T_{p,\alpha}}A(s)=1/2$. In any case, there
{\em does} exist some {\em maximal} $T_{p,\alpha}$ as above, and $A$
is uniquely defined as a function on $[0,T_{p,\alpha})$.

Our theorem can now be stated.
\begin{theorem}\label{thm:ODE_main}Consider a balls-in-bins process with feedback function $f(x)=x^p$, $0<p<1/2$ started from initial conditions $[t,\alpha]$. Let $\hat\alpha_t(s)$ be the fraction of balls in bin $1$ at time $\lceil s\,t\rceil$, and let $\elead$ be the event that bin $1$ eventually leads. It then holds that for all $K\in \R^+$ satisfying $K<T_{p,\alpha}$,
\begin{equation}\Prp{[t,\alpha]}{\forall s\in[0,K]\commaeq\;|\hat\alpha_t(s)-A_{p,\alpha}(s)|\leq W\, t^{-1/3}\mid \elead}\geq 1 - e^{-\ohmega{t^{1/3}}}\commaeq\; t\gg 1\periodeq\end{equation}Here, $W$ is a constant depending on $\alpha$ and $K$, but not on $t$.\end{theorem}

Notice that the two possibilities presented above --
$T_{p,\alpha}<+\infty$ or $T_{p,\alpha}=+\infty$ -- are both {\em a
priori} legitimate. In the former case, one would be able to show
that the random function $\hat\alpha(\cdot)$ conditioned on $\elead$
converges weakly to $A_{p,\alpha}$ in the space $D[0,+\infty)$
\cite{BillingsleyBook}. In the latter case, for all $\eps>0$, there
would be a value of $K=K_\eps<T_{p,\alpha}$ such that, with
probability tending to $1$ $\hat\alpha(K_\eps)>1/2 -\eps$. It would
be quite interesting to settle this matter: determining whether
$A(s)\to 1/2$ as $s$ converges to some finite $T$ should only
require a careful (but perhaps laborious) estimation of the RHS of
$g(\alpha)$ for $\alpha$ near $1/2$. The proof of \thmref{ODE_main}
can be found in \secref{ODE} below.

Before we proceed, let us briefly discuss our proof techniques. This
work employs the same fundamental tool as in the remaining papers in
this series \cite{MitzenmacherOS04,EuOnset,EuBrown}, as well as in
other references \cite{Khanin01,Spencer??} (according to
\cite{Khanin01}, the techique originated in Davis' work on
reinforced random walks \cite{Davis90}). We shall {\em embed} the
discrete-time process we are interested in into a continuous-time
process built from exponentially distributed random variables, so
that inter-arrival times at different bins are independent and have
an explicit distribution, which is very helpful in calculations. We
call this the {\em exponential embedding} of the process. Our main
conceptual contribution is to notice that the problems at hand lend
themselves to proof via the exponential embedding method.

The rest of the paper is organized as follows. After preliminaries
are discussed in \secref{prelim}, \secref{exp_embed} describes the
exponential embedding and its application to the eventual leadership
event. The next two sections prove the main theorems, and
\secref{final} discusses some open questions.

\section{Preliminaries}\label{sec:prelim}

\par {\em General notation.} Throughout the paper,
$\N=\{1,2,3,\dots\}$ is the set of non-negative integers, $\R^+ =
[0,+\infty)$ is the set of non-negative reals, and for any
$k\in\N\backslash\{0\}$ $[k]=\{1,\dots,k\}$. $\Ind{A}$ is the indicator function of a set (or event) $A$.\\

\par {\em Asymptotics.} We use the standard $O/o/\Omega/\Theta$ notation. The expressions ``$a_n\sim b_n$ as $n\to n_0$" and ``$a_n\ll b_n$ as $n\to n_0$" mean that $\lim_{n\to n_0}(a_n/b_n) = 1$ and $\lim_{n\to n_0}(a_n/b_n) = 0$, respectively.\\

\par {\em Balls-in-bins.} Formally, a feedback function is a map $f:\N\to (0,+\infty)$ with positive minimum. A balls-in-bins process with feedback function $f$ and $B\in\N$ bins is a discrete-time Markov chain $\{(I_1(m),\dots,I_B(m))\}_{m=0}^{+\infty}$ with state space $\N^B$ and transitions given as follows. For every time $m\geq 1$ there exists an index $i_m\in[B]$ such that $I_{m}(i_m)=I_{m-1}(i_m)+1$ and $I_{m}(i)=I_{m-1}(i)$ for
$i\in[B]\backslash\{i_m\}$. Moreover, the distribution of $i_m$ is
given by
$$\Pr{i_m=i\mid \{I_{m'}(j)\,:\, 0\leq m'<m, j\in[B]\}} =
\frac{f(I_{m-1}(i))}{\sum_{j=1}^B f(I_{m-1}(j))}.$$ We will usually
refer to the index $i_m\in[B]$ as {\em the bin that receives a ball
at time $m$}. For any $B$, if $E$ is an event of the process and
$u\in\N^B$, $\Prp{u}{E}$ is the probability of $E$ when the initial
conditions
are set to $u$. Finally, in the case $B=2$, it will be convenient to use the notation $[t,\alpha]$ ($t\in\N$, $0\leq \alpha\leq 1$ to denote the state $(\lceil\alpha t\rceil,t - \lceil\alpha t\rceil)$, i.e. there is a total of $t$ balls in the bins, and the fraction of balls in bin $1$ is (approximately) $\alpha$.\\

\par {\em Exponential random variables.} $X\eqdist \exp(\lambda)$
means that $X$ is a random variable with exponential distribution
with rate $\lambda>0$, meaning that $X\geq 0$ and
$$\Pr{X>t} = e^{-\lambda t}\;\; (t\geq 0).$$
The shorthand $\exp(\lambda)$ will also denote a generic random
variable with that distribution. Some elementary but extremely
useful properties of those random variables include
\begin{enumerate}\item {\em Lack of memory.} Let $X\eqdist \exp(\lambda)$
and $Z\geq 0$ be independent from $X$. The distribution of $X-Z$
conditioned on $X>Z$ is still equal to $\exp(\lambda)$. \item {\em
Minimum property.} Let $\{X_i\eqdist \exp(\lambda_i)\}_{i=1}^{m}$ be
independent. Then $$X_{\min}\equiv \min_{1\leq i\leq m}X_i\eqdist
\exp(\lambda_1+\lambda_2+\dots \lambda_m)$$ and for all $1\leq i\leq
m$
\begin{equation}\Pr{X_i = X_{\min}} = \frac{\lambda_i}{\lambda_1+\lambda_2+\dots
\lambda_m}\end{equation}  \item {\em Multiplication property.} If
$X\eqdist\exp(\lambda)$ and $\eta>0$ is a fixed number, $\eta X
\eqdist\exp(\lambda/\eta)$. \item {\em Moments and transforms.} If
$X\eqdist\exp(\lambda)$, $r\in\N$ and $t\in\R$,
\begin{eqnarray}\label{eq:moments}
\Ex{X^r} & = & \frac{r!}{\lambda^r}\commaeq \\
\Ex{e^{t X}} & = & \left\{\begin{array}{ll}\frac{1}{1 - \frac{t}{\lambda}} & (t<\lambda)\\
                    +\infty  & (t\geq
                    \lambda)\end{array}\right.\end{eqnarray}

\end{enumerate}

\section{The exponential embedding}\label{sec:exp_embed}

\subsection{Definition and key properties}

Let $f:\N\to(0,+\infty)$ be a feedback function, $B\in\N$ and
$(a_1,\dots,a_B)\in\N^B$. We define below a continuous-time process
with state space $(\N\cup\{+\infty\})^B$ and initial state
$(a_1,\dots,a_B)$ as follows. Consider a set $\{X(i,j)\,:\,i\in
[B],\, j\in \N\}$ of independent random variables, with
$X(i,j)\eqdist \exp(f(j))$ for all $(i,j)\in [B]\times \N$, and
define

\begin{equation}\label{eq:exp_contproc}N_i(t)\equiv \sup\left\{n\in
\N\,:\,\sum_{j=a_i}^{n-1} X(i,j) \leq t\right\}\;\;\;(i\in
[B],t\in\R^+ = [0,+\infty))\commaeq\end{equation} where by
definition $\sum_{j=i}^{k}(\dots)=0$ if $i>k$. Thus $N_i(0)=a_i$ for
each $i\in [B]$, and one could well have $N_i(T)=+\infty$ for some
finite time $T$ (indeed, that {\em will} happen for our cases of
interest); but in any case, the above defines a continuous-time
stochastic process, and in fact the $\{N_i(\cdot)\}_{i=1}^B$
processes are independent. Each one of this processes is said to
correspond to {\em bin} $i$, and each one of the times
$$X(i,a_i),X(i,a_{i})+X(i,a_{i}+1),X(i,a_{i})+X(i,a_{i}+1)+X(i,a_{i}+2),\dots$$
is said to be an {\em arrival time at bin $i$}. As in the
balls-in-bins process, we imagine that each arrival correspond to a
ball being placed in bin $i$.

In fact, we {\em claim} that this process is related as follows to
the balls-in-bins process with feedback function $f$, $B$ bins and
initial conditions $(a_1,\dots,a_B)$.

\begin{theorem}[Proven in \cite{Davis90,Khanin01,Spencer??,tese,EuOnset}]\label{thm:exp_embed}Let the $\{N_i(\cdot)\}_{i\in[B]}$ process be defined as above. One can order the arrival times of the $B$ bins in increasing order (up to their first accumulation point, if they do accumulate) so that $T_1<T_2<\dots$ is the resulting sequence. The distribution of
$$\{I_m = (N_1(T_m),N_2(T_m),\dots,N_B(T_m))\}_{m\in\N}$$
is the same as that of a balls-in-bins process with feedback
function $f$ and initial conditions
$(a_1,a_2,\dots,a_B)$.\end{theorem}

One can prove this result\footnote{The exact attribution of this
result is somewhat confusing. Ref. \cite{Khanin01} cites the work of
Davis \cite{Davis90} on reinforced random walks, where it is in turn
attributed to Rubin.} as follows. First, notice that the {\em first
arrival time $T_1$} is the minimum of $X(j,a_j)$, ($1\leq j\leq B$).
By the minimum property presented above, the probability that bin
$i$ is the one at which the arrival happens is like the first
arrival probability in the corresponding balls-in-bins process with
feedback:
\begin{equation}\label{eq:min}\Pr{X(i,a_i) = \min_{1\leq j\leq
B}X(j,a_j)} = \frac{f(a_i)}{\sum_{j=1}^B f(a_j)}.\end{equation} More
generally, let $t\in\R^+$ and condition on $(N_i(t))_{i=1}^B =
(b_i)_{i=1}^B \in\N^B$, with $b_i\geq a_i$ for each $i$ (in which
case the process has not blown up). This amounts to conditioning on
$$\forall i\in[B]\;\; \sum_{j=a_i}^{b_i-1}X(i,b_i)\leq t < \sum_{j=a_i}^{b_i}X(i,b_i).$$
From the lack of memory property of exponentials, one can deduce
that the first arrival after time $t$ at a given bin $i$ will happen
at a $\exp(f(b_j))$-distributed time, independently for different
bins. This almost takes us back to the situation of \eqnref{min},
with $b_i$ replacing $a_i$, and we can similarly deduce that bin $i$
gets the next ball with the desired probability,
$$\frac{f(b_i)}{\sum_{j=1}^Bf(b_j)}.$$

\subsection{On the eventual leadership event}\label{sec:exp_elead}

Before we move on to the main proofs, let us briefly discuss how the
event $\elead$ corresponding to eventual leadership by bin $1$ can
be expressed via the exponential embedding. We use the same notation
and random variables introduced above, and in particular we use the
embedded version of the balls-in-bins process defined above.
However, we restrict ourselves to the $B=2$ case with
\begin{equation}\label{eq:elead_condition}\sum_{j=1}^{+\infty}f(j)^{-2}<+\infty.\end{equation}
Notice that this condition implies we are either in the monopolistic
or in the eventual leadership regimes. Assume we start the process
from state $(x,y)\in\N^2$ with $x<y$ (i.e. bin $1$ has less balls
than bin $2$). The event $\elead$ is given by
$$\elead \equiv \{\exists m\geq 0 \forall M\geq m\,
I_m(1)>I_m(2)\}.$$ This can be restated as follows. For
$i\in\{1,2\}$, let $U^{(i)}_r$ be the first $m\in\N$ such that
$I_m(i)=r$, or set $U^{(i)}_r=+\infty$ if no such $m$ exists. Then

$$\elead \equiv \{\exists r\geq 0 \forall R\geq r\,
U^{(1)}_R<U^{(2)}_R\}.$$

This carries over to the {\em continuous-time process}, in which the
time it takes for bin $1$ to reach level $R$ is
$\sum_{j=x}^{R-1}X(1,j)$, and the analogous time for bin $2$ is
$\sum_{j=y}^{R-1}X(2,j)$. It is easy to show that
\begin{eqnarray}\elead &=& \left\{ \exists r\geq 0\, \forall R\geq r
\sum_{j=x}^{R-1}X(1,j)<\sum_{j=y}^{R-1}X(2,j)\right\} \\ &=&
\label{eq:elead_quase}\left\{ \exists r\geq 0\, \forall R\geq r
\sum_{j=y}^{R-1}(X(1,j)-X(2,j))
-\sum_{j=x}^{y-1}X(2,j)<0\right\}.\end{eqnarray}

The key now is to show that $\sum_{j=y}^{R-1}(X(1,j)-X(2,j))$
converges as $R\to +\infty$. Indeed, the random variables in the
sum,
$$X(i,j), i\in\{1,2\}, j\geq y$$
are independent, and each term in the sum has zero mean (since
$X(1,j)\eqdist X(2,j)\eqdist \exp(f(j))$) and variances that add up
to (cf. \eqnref{moments})
$$\sum_{j=y}^{R-1}\Var{X(1,j)-X(2,j)} = \sum_{j=y}^{R-1}\frac{2}{f(j)^2}\to \sum_{j=y}^{+\infty}\frac{2}{f(j)^2}<+\infty (\mbox{by \eqnref{elead_condition}}).$$
Kolmogorov's Three Series Theorem then implies that
$\sum_{j=y}^{+\infty}(X(1,j)-X(2,j))\in\R$ is a well-defined random
variable, as stated. Moreover, the event in \eqnref{elead_quase}
holds if and only if $\sum_{j=y}^{+\infty}(X(1,j)-X(2,j))
-\sum_{j=x}^{y-1}X(2,j)<0$, except for a null event, because
$\sum_{j=y}^{+\infty}(X(1,j)-X(2,j))$ and $\sum_{j=x}^{y-1}X(2,j)$
are independent (by the definition of the exponential embeddings)
and have no point masses in their distributions. It follows that

\begin{equation}\label{eq:elead_final}\Prp{(x,y)}{\elead} = \Prp{(x,y)}{\sum_{j=y}^{+\infty}(X(1,j)-X(2,j)) -
\sum_{j=x}^{y-1}X(2,j)<0}.\end{equation}

This equation is fundamental to our proofs.

\section{Escaping a very likely defeat}\label{sec:escape}

In this section we present the proof of \thmref{escape_main}. For
convenience, we have divided our argument into four parts. In
\subref{escape_method} we outline our proof method, which consists
of careful estimates of Laplace transforms. Such estimates are
carried out in \subref{escape_laplace} and
\subref{escape_asymptotic}. Those results are collected and applied
to the proof of the Theorem in \subref{escape_proof}.

\subsection{Our method of proof}\label{sub:escape_method}

As usual, our proof begins by writing down the event under
consideration in terms of the exponential embedding random
variables, using in this case \eqnref{elead_final} with
$(x,y)=[t,\alpha]$.
\begin{equation}\Prp{[t,\alpha]}{\elead} =\Pr{\sum_{j=\lceil\alpha t\rceil}^{t-\lceil\alpha
t\rceil-1}X(1,j) + \sum_{j=t-\lceil\alpha
t\rceil}^{+\infty}(X(1,j)-X(2,j))<0 }\periodeq\end{equation} Hence,
if we define the following independent random variables
\begin{eqnarray}A_t &\equiv & \sum_{j=\lceil\alpha t\rceil}^{t-\lceil\alpha
t\rceil-1}X(1,j)\commaeq\\
 \Delta_t &\equiv & \sum_{j=t-\lceil\alpha
t\rceil}^{+\infty}(X(1,j)-X(2,j))\commaeq \end{eqnarray} and let
$Z_t\equiv A_t+\Delta_t$, we deduce that
\begin{equation}\label{eq:escape_eleadgood2}\Prp{[t,\alpha]}{\elead} = \Pr{Z_t<0}\end{equation}
and that, for all $\lambda>0$,
\begin{equation}\label{eq:escape_eleadgood3}\Prp{[t,\alpha]}{\elead} \leq \Ex{\exp(-\lambda Z_t)} = \Ex{\exp(-\lambda A_t)} \Ex{\exp(-\lambda \Delta_t)}\periodeq\end{equation}

Thus the ``standard trick" of employing the Laplace transform
provides an upper bound on $\elead$. We now use a less standard
trick for lower bounding this probability in terms of the same
Laplace Transform. Our approach is essentially that of Spencer
\cite{Spencer??2}.

Let $\lambda>0$ and $\eta_1>\eta_2>0$ be given. Then
\begin{equation}\label{eq:escape_joelstrick}\Pr{Z_t<0} \geq \Pr{-\eta_1<Z_t<-\eta_2} \geq e^{-\lambda \eta_1}\,\Ex{e^{-\lambda\,Z_t}\Ind{\{-\eta_1<Z_t<-\eta_2\}}}\commaeq\end{equation}
since $-\eta_1<Z_t<\eta_2$ implies that $\lambda\eta_1 > -\lambda
Z_t$. Now let $\eps>0$ be fixed. Then, if $Z_t>-\eta_2$, then
$-\lambda Z_t <-(1-\eps)\lambda Z_t -\eta_2\eps$, and if
$Z_t<-\eta_1$, then $-\lambda Z_t < -\lambda (1 + \eps)Z_t -
\eps\lambda\eta_1$. Thus
\begin{multline}e^{-\lambda Z_t} < e^{-(1+\eps)\lambda\,Z_t - \lambda\eps\eta_1} +  e^{-(1-\eps)\lambda\,Z_t - \eps\lambda\eta_2} \\ \mbox{ on the complement of }\{-\eta_1<Z_t<-\eta_2\}\periodeq\end{multline}Hence
$$\Ex{e^{-\lambda\,Z_t}\Ind{Z_t\in (-\eta_1,-\eta_2)}} > \Ex{e^{-\lambda Z_t}} - e^{-\lambda\eps\eta_1}\,\Ex{e^{-(1+\eps)\lambda Z_t}} - \,e^{-\lambda\eps\eta_2}\Ex{e^{-(1-\eps)\lambda Z_t-\eps\lambda\eta}}$$
whenever the Laplace transforms above are finite. Plugging this last
inequality back into \eqnref{escape_joelstrick} yields the following
general lower bound.
\begin{multline}\label{eq:escape_joelstrick2}\Pr{Z_t<0} \\ \geq e^{-\lambda \eta_1}\Ex{e^{-\lambda\,Z_t}} \left(1 - \frac{e^{-\lambda\eps\eta_1}\,\Ex{e^{-(1+\eps)\lambda Z_t}} + \,e^{-\lambda\eps\eta_2}\Ex{e^{-(1-\eps)\lambda Z_t}}}{\Ex{e^{-\lambda Z_t}}}\right)\commaeq\end{multline}
How can one use the upper and lower bounds above? For the sake of
understanding what follows, let us indicate how inequalities
\eqnref{escape_eleadgood3} and \eqnref{escape_joelstrick2} are
typically employed. Assume that there is a choice of
$\lambda^*=\lambda^*_t$ that minimizes or nearly minimizes the
expression
\begin{equation}h_t(\lambda) \equiv \frac{1}{t}\ln\Ex{e^{-\lambda\,Z_t}}\end{equation}
Then one could hope that $h_t'(\lambda^*)\approx 0$,
$h_t''(\lambda^*)>0$, and that there would exist a constant $a$ {\em
not depending on $t$} such that for all $\delta>0$ small enough
\begin{equation}h_t((1\pm \delta)\lambda^*) \leq h_t(\lambda) + a\delta^2\periodeq\end{equation}
Thus our main expectation is that $h_t$ has an minimizer $\lambda^*$
and that it behaves like a ``nice" strictly convex function around
$\lambda^*$ in a way that does not depend on $t$. Now assume that
$\eps$ is small enough (but fixed) and we set
$\eta_1=\sqrt{\eps}t/\lambda$, $\eta_2=\sqrt{\eps}t/\lambda$ in
\eqnref{escape_joelstrick}, then
\begin{eqnarray}\frac{e^{-\lambda\eta_{1}\eps}\Ex{e^{-(1+\eps)\lambda^*\,Z_t}}}{\Ex{e^{-\lambda^*\,Z_t}}} &=& \exp\left\{[h_t((1+\eps)\lambda^*)-h_t(\lambda^*)-\eps\,]\,t\right\} \\
&\leq & \exp\{(a\eps^2 -\eps^{3/2}) t\} =
e^{-\ohmega{t}}\periodeq\end{eqnarray} and similarly
\begin{equation}\frac{e^{-\lambda\eta_{2}\eps}\Ex{e^{-(1-\eps)\lambda^*\,Z_t}}}{\Ex{e^{-\lambda^*\,Z_t}}} = e^{-\ohmega{t}}\periodeq\end{equation}
Thus in this case, \eqnref{escape_joelstrick} and
\eqnref{escape_eleadgood3} (with the choice of $\lambda=\lambda^*$)
would imply that
\begin{equation}(1-e^{-\ohmega{t}})e^{h_t(\lambda^*)t - \sqrt{\eps}t}\leq \Pr{Z_t<0}\leq e^{h_t(\lambda^*)t}\periodeq\end{equation}
This last expression would imply that
\begin{equation}h_t(\lambda^*) - \sqrt{\eps} -\liloh{1} \leq \frac{\ln\Pr{Z_t<0}}{t}\leq h_t(\lambda^*) +\liloh{1}\;\mbox{ for }t\gg 1\commaeq\end{equation}
for all small enough $\eps$, which shows that
\begin{equation}\lim_{t\to +\infty}\frac{\ln\Pr{Z_t<0}}{t} - h_t(\lambda^*)= 0\periodeq\end{equation}
The above exposition does not exactly correspond to our proof of
\thmref{escape_main}. However, the spirit of the two proofs is the
same. That is, we will show that the logarithms of our Laplace
transforms are ``strictly convex in the limit", and use that to
prove the desired result.

\subsection{Analysis of the Laplace
Transform}\label{sub:escape_laplace}

To apply the above method, we need to analyze the Laplace transform
of $Z_t=\Delta_t+A_t$. We start with $\Ex{\exp(-\lambda A_t)}$.
\begin{eqnarray}\Ex{\exp(-\lambda A_t)} &=& \prod_{j=\lceil\alpha t\rceil}^{t-\lceil\alpha t\rceil-1}\frac{1}{1 + \frac{\lambda}{j^p}}\\
& = & \exp\left\{\sum_{j=\lceil\alpha t\rceil}^{t - \lceil\alpha
t\rceil} \ln\left(\frac{1}{1 +
\frac{\lambda}{j^p}}\right)\right\}\end{eqnarray} With foresight, we
parameterize $\lambda = \lambda(\rho) = \rho (1-\alpha)^p t^p$, for
some $\rho>0$, and deduce that
\begin{eqnarray} &  &\Ex{\exp(-\lambda A_t)} = \\ &=& \exp\left\{\sum_{j=\lceil\alpha t\rceil}^{t - \lceil\alpha
t\rceil-1} \ln\left(\frac{1}{1 +
\frac{(1-\alpha)^p\rho}{(j/t)^p}}\right)\right\}\\
& = & \exp\left\{ t \times \left[\frac{1}{t}\,\sum_{j=\lceil\alpha
t\rceil}^{t - \lceil\alpha t\rceil-1} \ln\left(\frac{1}{1 +
\frac{(1-\alpha)^p\rho}{(j/t)^p}}\right)\right]\right\}\periodeq\end{eqnarray}
It is easy to see that the bracketed term is (close to) a Riemmann
sum. In fact, the function $u\mapsto
\ln(1/(1+(1-\alpha)^p\rho/u^p))$ is monotone increasing, so for any
$\lceil \alpha t\rceil \leq j\leq t - \lceil \alpha t\rceil $,
\begin{multline}0\leq \int_{\frac{j}{t}}^{\frac{j+1}{t}}\ln\left(\frac{1}{1 +
\frac{(1-\alpha)^p\rho}{v^p}}\right)\,dv
-\frac{1}{t}\ln\left(\frac{1}{1
+\frac{(1-\alpha)^p\rho}{(j/t)^p}}\right)\\ =
(1-\alpha)\int_{\frac{j}{(1-\alpha)t}}^{\frac{j+1}{(1-\alpha)t}}\ln\left(\frac{1}{1
+ \frac{\rho}{u^p}}\right)\,du - \frac{1}{t}\ln\left(\frac{1}{1
+\frac{(1-\alpha)^p\rho}{(j/t)^p}}\right) \\
\frac{1}{t}\left[\ln\left(\frac{1}{1
+\frac{(1-\alpha)^p\rho}{((j+1)/t)^p}}\right) - \ln\left(\frac{1}{1
+\frac{(1-\alpha)^p\rho}{(j/t)^p}}\right)
\right]\periodeq\end{multline} Summing over $j$ then yields
\begin{multline}0\leq
(1-\alpha)\int_{\frac{\lceil\alpha
t\rceil}{(1-\alpha)t}}^{\frac{t-\lceil\alpha
t\rceil}{(1-\alpha)t}}\ln\left(\frac{1}{1 +
\frac{\rho}{u^p}}\right)\,du - \frac{1}{t}\,\sum_{j=\lceil\alpha
t\rceil}^{t - \lceil\alpha t\rceil-1} \ln\left(\frac{1}{1 +
\frac{(1-\alpha)^p\rho}{(j/t)^p}}\right)\\ \leq
\frac{1}{t}\left[\ln\left(\frac{1}{1
+\frac{(1-\alpha)^p\rho}{((t-\lceil\alpha t\rceil)/t)^p}}\right) -
\ln\left(\frac{1}{1 +\frac{(1-\alpha)^p\rho}{(\lceil\alpha
t\rceil/t)^p}}\right) \right]\periodeq\end{multline} Thus we deduce
that for all $\delta>0$ there exist $C=C^{(1)}_\delta>0$,
$T_\delta\in\N$ such that if $\alpha>\delta$ and $t\geq
T^{(1)}_\delta$, then
\begin{multline}-\frac{C^{(1)}_\delta}{t}\leq \frac{1}{t}\,\sum_{j=\lceil\alpha
t\rceil}^{t - \lceil\alpha t\rceil-1} \ln\left(\frac{1}{1 +
\frac{(1-\alpha)^p\rho}{(j/t)^p}}\right) -
(1-\alpha)\int_{\frac{\alpha}{(1-\alpha)}}^{1}\ln\left(\frac{1}{1 +
\frac{\rho}{u^p}}\right)\,du \leq
\frac{C^{(1)}_\delta}{t}\periodeq\end{multline}

It follows that, for all $\alpha>\delta$ and $t\geq T^{(1)}_\delta$
\begin{multline}(1-\alpha)\int_{\frac{\alpha}{(1-\alpha)}}^{1}\ln\left(\frac{1}{1
+ \frac{\rho}{u^p}}\right) \leq\frac{\ln \Ex{\exp(-\lambda(\rho)
A_t)}}{t}  \\ \leq
(1-\alpha)\int_{\frac{\alpha}{(1-\alpha)}}^{1}\ln\left(\frac{1}{1 +
\frac{\rho}{u^p}}\right) +
\frac{C^{(1)}_\delta}{t}\periodeq\end{multline} We now consider the
Laplace transform of $\Delta_t$, with the same parametrization
$\lambda=\lambda(\rho)$ as above.
\begin{eqnarray}\Ex{\exp(-\lambda \Delta_t)} &=& \prod_{j=t-\lceil\alpha t\rceil}^{+\infty}\frac{1}{1-\frac{\lambda^2}{j^{2p}}} \\
&=& \exp\left\{\sum_{j=t-\lceil\alpha
t\rceil}^{+\infty}\ln\left(\frac{1}{1-\frac{(1-\alpha)^{2p}\rho^2}{(j/t)^{2p}}}\right)\right\}
\\
&=& \exp\left\{\left[\frac{1}{t}\sum_{j=t-\lceil\alpha
t\rceil}^{+\infty}\ln\left(\frac{1}{1-\frac{(1-\alpha)^{2p}\rho^2}{(j/t)^{2p}}}\right)\right]\,
t\right\}\periodeq\end{eqnarray} Notice that this Laplace transform
is infinite for $\lambda\geq t - \lceil\alpha t\rceil$, and we
therefore place the restriction $\rho\in(0,1)$ to ensure that does
not happen for all large enough $t$. As above, we have a something
close to a Riemmann sum between brackets. Indeed, since the map
$u\mapsto \ln(1/(1-(1-\alpha)^{2p}\rho^2/u^{2p}))$ is monotone
decreasing, for all $j\geq t - \lceil\alpha t\rceil$
\begin{multline}0\leq \frac{1}{t}\ln\left(\frac{1}{1-\frac{(1-\alpha)^{2p}\rho^2}{(j/t)^{2p}}}\right) -  \int_{\frac{j}{t}}^{\frac{j+1}{t}}\ln\left(\frac{1}{1-\frac{(1-\alpha)^{2p}\rho^2}{v^{2p}}}\right)\,dv \\
\leq
\frac{1}{t}\ln\left(\frac{1}{1-\frac{(1-\alpha)^{2p}\rho^2}{(j/t)^{2p}}}\right)
-(1-\alpha)\int_{\frac{j}{(1-\alpha)t}}^{\frac{j+1}{(1-\alpha)t}}\ln\left(\frac{1}{1-\frac{\rho^2}{u^{2p}}}\right)\,du
\\ \leq \left[\frac{1}{t}\ln\left(\frac{1}{1-\frac{(1-\alpha)^{2p}\rho^2}{(j/t)^{2p}}}\right) - \frac{1}{t}\ln\left(\frac{1}{1-\frac{(1-\alpha)^{2p}\rho^2}{((j+1)/t)^{2p}}}\right)\right]\periodeq\end{multline}
Summing over $j$, we conclude that
\begin{multline}0\leq \frac{1}{t}\sum_{j=t-\lceil\alpha t\rceil}^{+\infty}\ln\left(\frac{1}{1-\frac{(1-\alpha)^{2p}\rho^2}{(j/t)^{2p}}}\right) - (1-\alpha)\int_{\frac{t-\lceil\alpha t\rceil}{(1-\alpha)t}}^{+\infty}\ln\left(\frac{1}{1-\frac{\rho^2}{u^{2p}}}\right)\,du
\\ \leq \frac{1}{t}\ln\left(\frac{1}{1-\frac{(1-\alpha)^{2p}\rho^2}{((t-\lceil\alpha\rceil)/t)^{2p}}}\right)\periodeq\end{multline}
This implies that for each $\eta\in(0,1)$ there exist
$C^{(2)}_\eta>0$ and $T^{(2)}_\eta\in\N$ such that for all $t\geq
T^{(2)}_\eta$, if $0<\rho<1-\eta$
\begin{multline}-\frac{C^{(2)}_\eta}{t}\leq \frac{1}{t}\ln\Ex{\exp(-\lambda \Delta_t)} - (1-\alpha)\int_{\frac{t-\lceil\alpha t\rceil}{(1-\alpha)t}}^{+\infty}\ln\left(\frac{1}{1-\frac{\rho^2}{u^{2p}}}\right)\,du
\\ \leq \frac{C^{(2)}_\eta}{t}\periodeq\end{multline}
To conclude the section, let
\begin{eqnarray}\label{eq:escape_defgtrho}g_t(\rho,\alpha) &\equiv& \ln\Ex{\exp\{-\rho(1-\alpha)^pt^p[Z_t]\}} \\ \nonumber
& = & \sum_{j=\lceil\alpha t\rceil}^{t - \lceil\alpha t\rceil}
\ln\left(\frac{1}{1 + \frac{(1-\alpha)^p\rho}{(j/t)^p}}\right) \\
\nonumber & & + \sum_{j=t-\lceil\alpha
t\rceil}^{+\infty}\ln\left(\frac{1}{1-\frac{(1-\alpha)^{2p}\rho^2}{(j/t)^{2p}}}\right)\commaeq\;\alpha\in(0,1/2)\commaeq
\rho\in(0,1)\periodeq\end{eqnarray} Also define, for the same range
of $\alpha,\rho$,
\begin{eqnarray}F_p(\rho,\alpha) &\equiv& (1-\alpha)\int_{\frac{\alpha}{1-\alpha}}^{1} \ln\left(\frac{1}{1 +
\frac{\rho}{u^p}}\right)\,du \\ \nonumber &
&+(1-\alpha)\int_{1}^{+\infty}\ln\left(\frac{1}{1-\frac{\rho^2}{u^{2p}}}\right)\,du\periodeq\end{eqnarray}
From the above, we deduce that for any $\delta>0$, if we let
$C_{\delta,\eta}\equiv C^{(1)}_\delta+C^{(2)}_\eta$ then
$T_{\delta,\eta}\equiv T^{(1)}_\delta+T^{(2)}_\eta$, then
\begin{equation}\label{eq:escape_mainlaplace}\forall \alpha\in(\delta,1/2)\commaeq\rho\in(0,1)\commaeq t\geq T_\delta\commaeq \left|\frac{g_t(\rho,\alpha)}{t} - F_p(\rho,\alpha)\right| \leq \frac{C_{\delta,\eta}}{t}\periodeq\end{equation}

\subsection{The asymptotic form of the Laplace
transform}\label{sub:escape_asymptotic}

We now analyze the function $F_p(\rho,\alpha)$ introduced above, as
well as its minimum over $\rho$, which we will prove to precisely
the function in the statement of the Theorem.
\begin{equation}\label{eq:escape_fpgood}c_p(\alpha) = \inf_{\rho\in(0,1)}F_p(\rho,\alpha)\;\;(\alpha\in(0,1/2))\periodeq\end{equation}
We have the following formulae for all $\alpha\in(0,1/2)$
\begin{eqnarray}\lim_{\rho\searrow 0}F_p(\rho,\alpha) &=& 0\commaeq\\
\lim_{\rho\nearrow 1}F_p(\rho,\alpha) &=& +\infty \commaeq\\
\frac{1}{1-\alpha}\frac{\partial F_p}{\partial\rho}(\rho,\alpha) &=&
- \int_{\frac{\alpha}{1-\alpha}}^{1} \frac{1}{u^p + \rho}\,du
+\int_{1}^{+\infty}\frac{2\rho}{u^{2p}-\rho^2}\,du\commaeq\\
\frac{1}{1-\alpha}\frac{\partial^2 F_p}{\partial^2\rho}(\rho,\alpha)
&=& \int_{\frac{\alpha}{1-\alpha}}^{1} \frac{1}{(u^p + \rho)^2}\,du
+\int_{1}^{+\infty}\frac{2u^{2p} +
2\rho^2}{(u^{2p}-\rho^2)^2}\,du\periodeq\end{eqnarray} Notice, then,
that
\begin{multline}\forall\alpha>0 \inf_{\rho\in(0,1)} \frac{\partial^2 F_p}{\partial^2\rho}(\rho,\alpha) \geq \inf_{\rho\in(0,1)}(1-\alpha)\int_{\frac{\alpha}{1-\alpha}}^{1} \frac{1}{(u^p + \rho)^2}\,du \\
=\inf_{\rho\in(0,1)}(1-\alpha)\int_{\frac{\alpha}{1-\alpha}}^{1}
\frac{1}{(u^p + 1)^2}\,du \equiv a_{\alpha}>0\commaeq\end{multline}
hence $F_p(\cdot,\alpha)$ is a strictly convex function of $\rho$,
for any $\alpha\in(0,1/2)$. Moreover,
\begin{equation}\lim_{\rho\searrow 0} \frac{\partial F_p}{\partial \rho}(\rho,\alpha) = -
(1-\alpha)\int_{\frac{\alpha}{1-\alpha}}^{1}
\frac{1}{u^p}\,du<0\commaeq\end{equation} The two last assertions
prove that for any fixed $\alpha$:
\begin{enumerate}
\item $F_p(\cdot,\alpha)$ is a strictly convex function of $\rho$;
\item $F_p(\rho,\alpha)<0$, $\partial_\rho F_p(\rho,\alpha)<0$ for
all small enough $\rho$; \item $F_p(\cdot,\alpha)$ thus has an
unique minimum over $(0,1)$, and this minimum is achieved at the
unique value $\rho^*=\rho^*(\alpha)$ such that
\begin{equation}\label{eq:escape_crhostar}\frac{\partial F_p}{\partial \rho}(\rho^*(\alpha),\alpha) = 0\semicoloneq\end{equation}
\item by the definition of $c_p(\alpha)$ and the above items,
\begin{equation}\label{eq:escape_fpgood2}c_p(\alpha) = \inf_{\rho\in(0,1)}F_p(\rho,\alpha) = F_p(\rho^*(\alpha),\alpha)<0\semicoloneq\end{equation}
\item by strict convexity (and using the definition of
$a_\alpha>0$ above) and the fact that $\partial_\rho
F_p(\rho^*(\alpha),\alpha)=0$, then there exists a value
$b=b(\alpha)$ depending continuously on $\alpha$ such that, for all
$\eps>0$ small enough,
\begin{equation}F_p((1\pm \eps)\rho^*(\alpha),\alpha)\leq F_p(\rho^*(\alpha),\alpha) +  b_\alpha\eps^2 = c_p(\alpha) + b_\alpha\eps^2\semicoloneq\end{equation}
\item in fact, we can strengthen the previous item and say that
for all $\delta\in(0,1/2)$ there exists $B_\delta,\eps_\delta>0$
such that
\begin{equation}\label{eq:escape_convexity}\forall \alpha\in(\delta,1/2-\delta)\commaeq \eps\in(0,\eps_\delta)\commaeq F_p((1\pm \eps)\rho^*(\alpha),\alpha)\leq c_p(\alpha) + B_\delta\eps^2\semicoloneq\end{equation}
\end{enumerate}
We now prove that $c_p$ is a smooth function of $\alpha$ with a
positive derivative. To prove smoothness, we only need to show that
$\rho^*$ is a smooth function of $\alpha$, since $c_p(\alpha)$ is
given by the formula in \eqnref{escape_crhostar}. But recall that we
have shown that $\rho^*$ is uniquely defined by the equation
\begin{equation}G(\rho^*(\alpha),\alpha) \equiv \frac{\partial F_p}{\partial \alpha}(\rho^*(\alpha),\alpha) = 0\commaeq\end{equation}
and we know that
\begin{equation}\frac{\partial G}{\partial\rho}(\rho,\alpha) =\frac{\partial^2 F_p}{\partial^2 \rho}(\rho,\alpha) \geq a_\alpha>0\periodeq\end{equation}
Hence the Implicit Function Theorem applies \cite{RudinBook}, and
implies that $\rho^*$ is indeed a smooth function of $\alpha$.

To prove that $c_p'(\alpha)>0$, we first differentiate
$F_p(\rho,\alpha)$ with respect to $\alpha$.
\begin{equation}\label{eq:escape_partiala}\frac{\partial F_p}{\partial\alpha}(\rho,\alpha) = (1-2\alpha)\log\left(1 + \frac{\rho}{\left(\frac{\alpha}{1-\alpha}\right)^p}\right) - \frac{1}{1-\alpha} F_p(\rho,\alpha)\periodeq\end{equation}
Now notice that, by the chain rule,
\begin{eqnarray}\label{eq:escape_partiala2}c_p'(\alpha) &=& \frac{d}{d\alpha}(F_p(\rho^*(\alpha),\alpha)) \\
&=& \frac{\partial F_p}{\partial\rho}(\rho^*(\alpha),\alpha) \,
(\rho^*)'(\alpha) + \frac{\partial
F_p}{\partial\rho}(\rho^*(\alpha),\alpha) \\
(\mbox{via }[\partial_\rho F_p](\rho^*,\alpha)=0)&=& \frac{\partial
F_p}{\partial\rho}(\rho^*(\alpha),\alpha)\\
(\mbox{by \eqnref{escape_partiala}}) & = &(1-2\alpha)\log\left(1 +
\frac{\rho}{\left(\frac{\alpha}{1-\alpha}\right)^p}\right) -
\frac{F_p(\rho,\alpha)}{1-\alpha}\\
(\mbox{by item $4.$ above}) & = & (1-2\alpha)\log\left(1 +
\frac{\rho}{\left(\frac{\alpha}{1-\alpha}\right)^p}\right)- \frac{c_p(\alpha)}{1-\alpha} \\
(\mbox{omitting a $>0$ term}) &\geq & - \frac{c_p(\alpha)}{1-\alpha} \\
(\mbox{since $c_p<0$}) &>& 0\periodeq\end{eqnarray}

Finally, we show that
\begin{equation}\label{eq:escape_rhodecreasing}\rho^*(\alpha) \mbox{ is a non-increasing function of $\alpha$.}\end{equation} This is important because it implies that, for all $\alpha>\delta>0$ and all small enough $\eps$,
$$(1+\eps)\rho^*(\alpha)\leq (1+\eps)\rho^*(\delta)\leq \eta<1$$
for some $\eta=\eta_\delta$ depending on $\delta$ only. In
conjunction with \eqnref{escape_mainlaplace}, this will imply that
\begin{multline}\label{eq:escape_mainlaplace2}\forall \alpha\in(\delta,1/2)\commaeq\rho\in(0,1)\commaeq t\geq T_\delta\commaeq \eps\in[0,\eps_\delta]\\ \left|\frac{g_t((1\pm \eps)\rho^*(\alpha),\alpha)}{t} - F_p((1\pm\eps)\rho^*(\alpha),\alpha)\right| \leq \frac{C_{\delta}}{t}\commaeq\end{multline}
where $C_\delta\equiv C_{\delta,\eta_\delta}$ depends on $\delta$
only.

To prove \eqnref{escape_rhodecreasing}, we notice that
\begin{eqnarray} & \frac{\partial^2
F_p}{\partial\alpha\partial\rho}(\rho^*(\alpha),\alpha)\\ \nonumber
= & \frac{\partial}{\partial\rho}\left[(1-2\alpha)\log\left(1 +
\frac{\rho}{\left(\frac{\alpha}{1-\alpha}\right)^p}\right) -
\frac{1}{1-\alpha}F_p(\rho,\alpha)\right]_{\rho=\rho^*(\alpha)}\\
\nonumber =&(1-2\alpha) \frac{1}{\frac{\alpha^p}{(1-\alpha)^p} +
\rho^*(\alpha)} - \frac{1}{1-\alpha}\frac{\partial
F_p}{\partial\rho}(\rho^*(\alpha),\alpha) \\
\label{eq:escape_bonitinha}=& (1-2\alpha)
\frac{1}{\frac{\alpha^p}{(1-\alpha)^p} +
\rho^*(\alpha)}>0\commaeq\end{eqnarray} where
\eqnref{escape_bonitinha} follows from
$$[\partial_\rho F_p](\rho^*(\alpha),\alpha)=0\periodeq$$
Hence, if $\beta>0$ is small enough
\begin{equation}\frac{\partial F_p}{\partial\rho}(\rho^*(\alpha),\alpha+\beta)>
\frac{\partial
F_p}{\partial\rho}(\rho^*(\alpha),\alpha)=0\periodeq\end{equation}
But by the strict convexity of $F_p(\cdot,\alpha+\beta)$,
$$\frac{\partial F_p}{\partial\rho}(\rho,\alpha+\beta)<\frac{\partial F_p}{\partial\rho}(\rho^*(\alpha+\beta),\alpha+\beta)$$
for all $\rho<\rho^*(\alpha+\beta)$. Hence
$\rho^*(\alpha+\beta)<\rho^*(\alpha)$ whenever $\beta$ is small
enough. This finishes the proof.

\subsection{Proof of \thmref{escape_main}}\label{sub:escape_proof}

We now have all the tools necessary to prove \thmref{escape_main}.

\begin{proof}[of \thmref{escape_main}] Let us now apply the upper and lower bounds
\eqnref{escape_eleadgood2} and \eqnref{escape_joelstrick2} presented
above. We will assume that $\delta\leq\alpha\leq 1/2-\delta$ for
some constant $\delta>0$, and prove bounds on
$\Prp{[t,\alpha]}{\elead}$ that are uniform on that range of
$\alpha$.

In the current setting, $Z_t=\Delta_t+A_t$ and we have defined
$$g_t(\rho,\alpha)= \ln\Ex{e^{-\lambda Z_t}}|_{\lambda = \rho[(1-\alpha)t]^p}\commaeq$$ hence for all fixed $\rho\in(0,1)$, $\alpha\in(\delta,1/2-\delta)$ and $t\geq T_\delta$ (cf. \eqnref{escape_mainlaplace})
\begin{equation}\Prp{[t,\alpha]}{\elead} = \Pr{Z_t<0}\leq \exp(g_t(\rho,\alpha))\periodeq\end{equation}
In particular, setting $\rho=\rho^*(\alpha)$ and $t$, $\alpha$ as
above, we can use \eqnref{escape_mainlaplace2} to bound
\begin{equation}\label{eq:escape_finalupper}\ln\Prp{[t,\alpha]}{\elead}\leq  \exp(c_p(\alpha)\,t + C_\delta)\periodeq\end{equation}
The above upper bound can be matched via the lower bound method in
\eqnref{escape_joelstrick2}. Let $\eps>0$. One can set $\lambda =
\rho^*(\alpha)[(1-\alpha)t]^p$, $\eta_1=\sqrt{\eps}t/\lambda$ and
$\eta_2=\sqrt{\eps}t/2\lambda$ in \eqnref{escape_joelstrick2} to
deduce
\begin{multline}\label{eq:escape_almostfinallower}\Prp{[t,\alpha]}{\elead} = \Pr{Z_t<0} \\
\geq e^{g_t(\rho^*,\alpha) -\sqrt{\eps}t}\left(1 -
\frac{e^{g_t((1+\eps)\rho^*,\alpha)-\eps^{3/2}t}}{e^{g_t(\rho^*)}}
-\frac{e^{g_t((1-\eps)\rho^*,\alpha)-\eps^{3/2}t/2}}{e^{g_t(\rho^*)}}\right)\periodeq\end{multline}
Now notice that and any $0<\eps<\eps_\delta$ (cf.
\eqnref{escape_convexity} and \eqnref{escape_mainlaplace2}), one has
that
\begin{eqnarray}g_t(\rho^*,\alpha) &=& c_p(\alpha)\,t \pm C_\delta\commaeq\\
g_t((1+\eps)\rho^*,\alpha) &=& F_p((1+\eps)\rho^*,\alpha)\, t \pm C_\delta \\
\nonumber (\mbox{by \eqnref{escape_convexity} and
\eqnref{escape_mainlaplace2}})&\leq&
c_p(\alpha)\,t + B_\alpha \eps^2\, t + C_\delta  \\
\nonumber (\mbox{for small enough $\eps$})&\leq& c_p(\alpha)\,t+ \frac{\eps^{3/2}}{4}\, t + C_\delta\commaeq \\
g_t((1-\eps)\rho^*,\alpha) &=& F_p((1-\eps)\rho^*,\alpha)\,t + C_\delta\\
\nonumber (\mbox{by \eqnref{escape_convexity} and
\eqnref{escape_mainlaplace2}})&\leq&
c_p(\alpha)\,t + B_\alpha \eps^2\, t + C_\delta  \\
\nonumber (\mbox{for small enough $\eps$}) &\leq& c_p(\alpha)\,t+
\frac{\eps^{3/2}}{4}\, t + C_\delta\periodeq\end{eqnarray}
Substitution back into \eqnref{escape_almostfinallower} yields
\begin{equation}\label{eq:escape_quasequase}\Prp{[t,\alpha]}{\elead} \geq (1-2e^{C_\delta-\eps^{3/2}t/4})\exp\{c_p(\alpha) t - C_\delta - \sqrt{\eps}t\}\commaeq\;t\geq T_\delta \\
\commaeq\end{equation} for any small enough $\eps$ and any $t\geq
T_\delta$. In particular, if we set $\eps \equiv [(C_\delta + \ln
4)(4t)]^{-2/3}$, then $\eps\searrow 0$ as $t\to +\infty$, so that
for all large enough $t$ the above formulae apply and
\begin{equation}\label{eq:escape_finallower}\Prp{[t,\alpha]}{\elead} \geq \frac{\exp(c_p(\alpha) t - C'_\delta\,t^{2/3})}{2}\commaeq\end{equation}
for some constant $C'_\delta\geq C_\delta$. Redefining $C_\delta$
and $T_\delta$ if necessary, we can then conclude (using
\eqnref{escape_finalupper} and \eqnref{escape_finallower}) that
\begin{equation}\label{eq:escape_final}\forall \alpha\in(\delta,1/2-\delta)\commaeq\forall t\geq T_\delta\commaeq\;\; c_p(\alpha)\, t - C_\delta\,t^{2/3} \leq \ln\Prp{[t,\alpha]}{\elead} \leq c_p(\alpha)\, t + C_\delta\periodeq\end{equation}
Since we have already shown that $c_p$ is smooth and
monotone-increasing in $\alpha$ (cf. \subref{escape_asymptotic}),
the Theorem follows.\end{proof}

\section{The most likely escape path}\label{sec:ODE}
This section is dedicated to \thmref{ODE_main}. After some
preliminaries are considered in \subref{ODE_prelim}, we then
estimate (in \subref{ODE_single}) the transition probabilities of
the balls-in-bins process conditioned on $\elead$. Those estimates
are used to show in \subref{ODE_central} that for short enough
times, the conditioned process evolves from a state $[t,\alpha]$ to
a state $\approx[(1+\eta)t,\alpha + g_p(\alpha)\eta]$, thus staying
close to the tangent of the ODE. The final steps of the proof are
presented in \subref{ODE_proof}, and a Lemma used in
\subref{ODE_central} is proven in \subref{ODE_Chernoff}.

\subsection{Preliminaries}\label{sub:ODE_prelim}

According to \thmref{escape_main}, the map
\begin{equation}\funcdef{c_p}{\left(0,\frac{1}{2}\right)}{\R^-}{\alpha}{\lim_{t\to +\infty}\frac{\ln\Prp{[t,\alpha]}{\elead}}{t}}\end{equation}
is infinitely differentiable. In particular, this means that, for
all $\delta\in(0,1/2)$, the suprema
\begin{equation}D^{(r)}_{\delta} \equiv \sup_{\delta\leq\alpha\leq \frac{1}{2}-\delta}\left|\frac{1}{r!}\frac{d^rc_p}{d\alpha^r}(\alpha)\right| \; (r\in\N\cup\{0\})\end{equation}
are all finite. Moreover, $c_p'(\alpha)>0$ on $(0,1/2)$, and
\begin{equation}\label{eq:ODE_minder}d^{(1)}_\delta \equiv \min_{\delta\leq \alpha\leq \frac{1}{2}-\delta}c_p'(\alpha) > 0\periodeq\end{equation}
\thmref{escape_main} also tells us that, for $\delta$ as above,
there exist $T_\delta\in\N$ and $C_\delta\in\R^+$ such that for all
$t\geq T_\delta$ and all $\alpha\in(\delta,1/2-\delta)$,
\begin{equation}- C_\delta t^{2/3}\leq \ln\Prp{[t,\alpha]}{\elead}- c_p(\alpha)\, t\leq C_\delta\periodeq\end{equation}
Therefore, if $\alpha,\alpha'\in(\delta,1/2-\delta)$ and $t,t'\geq
T_\delta$
\begin{equation}e^{c_p(\alpha') t' - c_p(\alpha)\, t - C_\delta (1+t^{2/3})}\leq \frac{\Prp{[t',\alpha']}{\elead}}{\Prp{[t,\alpha]}{\elead}} \leq e^{c_p(\alpha') t' - c_p(\alpha)\, t - C_\delta (1+(t')^{2/3})}\end{equation}
Moreover, if we also have that $\alpha-\alpha'=\eps$, $t'-t=\eta t$
for $\eps,\eta<1/4$ (say), then
\begin{multline}\label{eq:ODE_goodbound}\left|\frac{1}{t}\ln\frac{\Prp{[t',\alpha']}{\elead}}{\Prp{[t,\alpha]}{\elead}} - \eta\,c_p(\alpha) + \eps\,c_p'(\alpha)(1+\eta) \right| \\ \leq D^{(2)}_\delta\,\eps^2 + D^{(1)}_\delta \eps\eta + C_\delta\left(\frac{1}{t} + \frac{1}{t^{1/3}}\right)\periodeq\end{multline}
This last equation will be repeatedly used in what follows.

\subsection{Transitions conditioned on
$\elead$}\label{sub:ODE_single}

Define
\begin{equation}[t,\alpha]\mapsto [t^*,\alpha^*]\;\;\alpha,\alpha^*\in(0,1)\commaeq\;t\commaeq t^*\in\N,t<t^*\commaeq\end{equation}
to be the event that the initial state of the process is
$[t,\alpha]$, and that at time $t^*-t$ the state of the process is
$[t^*,\alpha^*]$. The goal of this section is to estimate the
probability of the transition
\begin{equation}\Prp{[t,\alpha]}{[t,\alpha]\mapsto [t^*,\alpha^*] \mid \elead}\end{equation}
for the case when $t\leq t^*\leq (1 +\eta)t$ for some $\eta>0$ and
$t$ is large, and we assume (for simplicity) that $\alpha t$,
$\alpha^*t^*$ are integers. We will require that
$\alpha\pm\eta\in(\delta,1/2-\delta)$ for some $0<\delta<1/2$, so
that for all $c\in(\alpha-\eta,\alpha+\eta)$ the bounds on
$\Prp{[t,c]}{\elead}$ coming from the previous section apply with
the same value of $\delta$ (and thus the same $C_\delta$,
$T_\delta$). Notice that we can assume that
\begin{equation}\label{eq:ODE_step1}\delta\leq \alpha - \eta\leq
\frac{\alpha}{1+\eta} \leq \alpha^* \leq \alpha + \eta\leq
1/2-\delta\commaeq\end{equation} otherwise the given probability is
$0$. One has that
\begin{eqnarray}\label{eq:ODE_rewrite} & & \Prp{[t,\alpha]}{[t,\alpha]\mapsto [t^*,\alpha^*] \mid \elead} \\ \nonumber &=& \frac{\Prp{[t,\alpha]}{\elead\mid [t,\alpha]\mapsto [t^*,\alpha^*]}\,\Prp{[t,\alpha]}{[t,\alpha]\mapsto [t^*,\alpha^*]}}{\Prp{[t,\alpha]}{\elead}}\\ \nonumber
&=&
\frac{\Prp{[t^*,\alpha^*]}{\elead}\,\Prp{[t,\alpha]}{[t,\alpha]\mapsto
[t^*,\alpha^*]}}{\Prp{[t,\alpha]}{\elead}}\commaeq\end{eqnarray}
where the first line is Bayes' Rule, and the second follows from the
Markov Property of the balls-in-bins process. Using the bounds in
\eqnref{ODE_step1} and \eqnref{ODE_goodbound},
\begin{equation}\frac{\Prp{[t^*,\alpha^*]}{\elead}}{\Prp{[t,\alpha]}{\elead}}= \exp(c_p(\alpha)(t^*-t) + c_p'(\alpha)(\alpha^*-\alpha)t^* \pm (D^{(2)}_\delta+ D^{(1)}_\delta) \eta^2\,t)\periodeq\end{equation}
It will be convenient to have the above equation in a slightly
different form,
\begin{equation}\label{eq:ODE_cansei}\frac{\Prp{[t^*,\alpha^*]}{\elead}}{\Prp{[t,\alpha]}{\elead}}= \exp(c_p(\alpha)(t^*-t) - \alpha c_p'(\alpha)(t^*-t) \pm (D^{(2)}_\delta+ D^{(1)}_\delta)\eta^2\,t) (e^{c_p'(\alpha)})^{\alpha^*t^* - \alpha t}\periodeq\end{equation}
As for
\begin{equation}\Prp{[t,\alpha]}{[t,\alpha]\mapsto[t^*,\alpha^*]}\commaeq\end{equation}notice that there $t^*-t = \eta t$ and that $\alpha^* t^* - \alpha t$, hence there exist
$$\binom{t^*-t}{\alpha^* t^* - \alpha t}$$
ways of moving from state $[t,\alpha]$ to state $[t^*,\alpha^*]$.
For each one of those ways, each step in which a ball is added to
bin $1$ has probability
$$\frac{(\alpha t + a)^p}{(\alpha t + a)^p + ((1-\alpha)t + b)^p}$$
of occurring (for some $0\leq a,b\leq t^*-t\leq \eta t$), whereas
steps in which a ball is added to bin $2$ have probability
$$\frac{((1-\alpha) t + b)^p}{(\alpha t + a)^p + ((1-\alpha)t + b)^p}\commaeq$$
for $a,b$ as above. There exist absolute constants $R_\delta$ and
$\eta_\delta$ only depending on $\delta$ such that, if
$0<\eta<\eta_\delta$:
$$e^{-R_\delta\,\eta}\,\frac{\alpha^p}{\alpha^p + (1-\alpha)^p}\leq \frac{(\alpha t + a)^p}{(\alpha t + a)^p + ((1-\alpha)t + b)^p}\leq e^{R_\delta\,\eta}\,\frac{\alpha^p}{\alpha^p + (1-\alpha)^p}$$
and
$$e^{-R_\delta\,\eta}\frac{(1-\alpha)^p}{\alpha^p + (1-\alpha)^p} \leq\frac{((1-\alpha) t + b)^p}{(\alpha t + a)^p + ((1-\alpha)t + b)^p}\leq e^{R_\delta\,\eta}\,\frac{(1-\alpha)^p}{\alpha^p + (1-\alpha)^p}\commaeq$$
and therefore, any path moving connecting state $[t,\alpha]$ to
state $[t^*,\alpha^*]$ has probability
\begin{equation}=e^{\pm R_\delta\,\eta (t^*- t)}\left(\frac{\alpha^p}{\alpha^p + (1-\alpha)^p} \right)^{\alpha^* t^* - \alpha t} \left(1 - \frac{\alpha^p}{\alpha^p + (1-\alpha)^p}\right)^{t^* - t - (\alpha^* t^* - \alpha t)} \commaeq \end{equation}
(any such path must have $\alpha^*t^*-\alpha t$ steps bin $1$
receives a ball, out of a total of $t^*-t$ steps). Letting
$$\rho(\alpha)\equiv\frac{\alpha^p}{\alpha^p + (1-\alpha)^p}\commaeq$$
we conclude that
\begin{multline}\Prp{[t,\alpha]}{[t,\alpha]\mapsto[t^*,\alpha^*]} =  e^{\pm R_\delta\,\eta^2 t}\\ \times\binom{t^*-t}{\alpha^*t^* - \alpha t}\rho(\alpha)^{\alpha^* t^* - \alpha t} \left(1 - \rho(\alpha)\right)^{(t^*-t) - (\alpha^* t^* - \alpha t)}\commaeq\end{multline}
and that (cf. \eqnref{ODE_cansei})
\begin{multline}\label{eq:ODE_javejo}\Prp{[t,\alpha]}{[t,\alpha]\mapsto[t^*,\alpha^*]\mid\elead} \\ =  \exp\{c_p(\alpha)(t^*-t) - \alpha c_p'(\alpha)(t^*-t)  \pm (R_\delta + D^{(2)}_\delta+ D^{(1)}_\delta)\,\eta^2 t\}\\ \times\binom{t^*-t}{\alpha^*t^* - \alpha t}[\rho(\alpha)e^{c_p'(\alpha)}]^{\alpha^* t^* - \alpha t} \left(1 - \rho(\alpha)\right)^{(t^*-t) - (\alpha^* t^* - \alpha t)}\periodeq\end{multline}
for all $t\geq T_\delta$, $0<\eta<\eta_\delta$, $t\leq t^*\leq
(1+\eta)t$ and $\alpha^*$ as above.

\subsection{The most likely transitions}\label{sub:ODE_central}

We continue with the same setup as above, and state a useful lemma
that we prove in \subref{ODE_Chernoff}.
\begin{lemma}\label{lem:ODE_Chernoff}Let $a>0$, $0<\rho<1$ and an integer $m\in\N$, $m\geq 2$ be given. Define, for $n\in [m]\cup\{0\}$,
$$b(n)\equiv\binom{m}{n}\rho^na^n(1-\rho)^{m-n}$$
Then\begin{enumerate} \item the sequence $\{b(n)\mid 0\leq n\leq
m\}$ is unimodal; \item $\max_{0\leq n\leq m}b(n)$ is achieved at
$n_0 = \left\lceil \frac{\rho am + (1-\rho)}{\rho
a+(1-\rho)}\right\rceil $; \item for all $K>0$,
$\sum_{|n-n_0|>K\sqrt{m}}b(n) \leq m\,b(n_0)
e^{-\frac{K^2}{2(1-n_0/m)}}$.\end{enumerate}\end{lemma} To apply
this lemma, we give new names to familiar quantities.
\begin{eqnarray} m=m(t^*,t) &=& t^*-t \leq \eta t\\
n=n(t,t^*\alpha,\alpha^*) &= & \alpha^*t^* - \alpha t\\
a=a(\alpha) &=& e^{c_p'(\alpha)}\\
\rho &=& \rho(\alpha)\\
n_0=n_0(\alpha,t^*,t) &=& \left\lceil\frac{\rho a m +
(1-\rho)}{(\rho a + 1-\rho)} \right\rceil\periodeq\end{eqnarray}
With those definitions,
$$\frac{n_0}{m} \leq \frac{\rho a m + (1-\rho)}{(\rho a + 1-\rho)m} + \frac{1}{m}\leq 1 - \frac{(1-\rho)(m-1)}{(\rho a + 1-\rho)m}\periodeq$$
For $\alpha\in(\delta.1/2-\delta)$, $a\leq e^{d^{(1)}_\delta}$.
Moreover, $\rho(\alpha)$ is a continuous function of $\alpha$ that
is between $\rho(\delta)>0$ and $1/2$ for $\alpha$ as above. Thus
there exists a constant $U=U_\delta\in(0,1)$ such that
$$\frac{n_0}{m} \leq 1 - U_\delta$$
and therefore the Lemma implies that, for $K>0$
\begin{equation}\sum_{|n-n_0|>K\sqrt{m}}b(n) \leq mb(n_0)e^{-\frac{K^2}{2U_\delta}}\periodeq\end{equation}
Now notice that, in the present case, \eqnref{ODE_javejo} implies
that
\begin{equation}b(n) = \frac{\Prp{[t,\alpha]}{[t,\alpha]\mapsto[t^*,\alpha]\mid \elead}}{\exp\{c_p(\alpha)(t^*-t) - \alpha c_p'(\alpha)(t^*-t)  \pm (R_\delta + D^{(2)}_\delta+ D^{(1)}_\delta)\,\eta^2 t\}}\periodeq\end{equation}
It follows that for all $K>0$, $\alpha\in(\delta,1/2-\delta)$,
$t\geq T_\delta$, $0<\eta<\eta_\delta$, $t\leq t^*\leq (1+\eta)t$
and $\alpha^*$ as above
\begin{multline}\label{eq:ODE_almostfinalbound}\Prp{[t,\alpha]}{[t,\alpha]\mapsto [t^*,\alpha]\mbox{ with }|n(t,t^*,\alpha,\alpha^*)-n_0(\alpha,t^*,t)|>K\sqrt{t^*-t}\mid\elead} \\
\leq (t^*-t)\exp\left\{2(R_\delta + D^{(2)}_\delta+
D^{(1)}_\delta)\,\eta^2 t - \frac{K^2}{2U_\delta}
\right\}\periodeq\end{multline} To use this formula, we will assume
that $t^*-t = \eta t$ (i.e. equality instead of the above
inequality), and then set $K = \sqrt{\eta}t$. Then
\begin{multline}\label{eq:ODE_almostfinalbound2}\Prp{[t,\alpha]}{[t,\alpha]\mapsto [t^*,\alpha]\mbox{ with }|n(t,t^*,\alpha,\alpha^*)-n_0(\alpha,t^*,t)|>\eta^{3/2} t \mid\elead} \\
\leq \eta t\,\exp\left\{2(R_\delta + D^{(2)}_\delta+
D^{(1)}_\delta)\,\eta^2 t - \frac{\eta\,t}{2U_\delta}
\right\}\commaeq\end{multline} and (by making $\eta_\delta$ smaller
if necessary) we can ensure that there exists $V_\delta>0$ such
that, with $0<\eta<\eta_\delta$,
\begin{multline}\label{eq:ODE_almostfinalbound3}\Prp{[t,\alpha]}{[t,\alpha]\mapsto [t^*,\alpha]\mbox{ with }|n(t,t^*,\alpha,\alpha^*)-n_0(\alpha,t^*,t)|>\eta^{3/2} t \mid\elead} \\
\leq e^{-V_\delta\, \eta\, t}\periodeq\end{multline} To conclude
this part, we look at how $\alpha^*$ behaves when
$$|n(t,t^*,\alpha,\alpha^*)-n_0(\alpha,t^*,t)|\leq \eta^{3/2} t\periodeq$$
In that case,
$$\alpha^* = \frac{\alpha t + n}{t^*} = \frac{\alpha}{1+\eta} + \frac{1}{1+\eta}\frac{n}{t}\periodeq$$
As defined above,
$$\frac{n_0}{t} = \frac{1}{t}\left\lceil\frac{\rho a \eta t + (1-\rho)}{(\rho a + 1-\rho)}\right\rceil = \frac{\rho a\eta }{\rho a + 1-\rho} \pm \frac{1-\rho}{(\rho a + 1-\rho)t}\periodeq$$
Hence, by making $\eta_\delta$ even smaller if necessary, we can
guarantee that
\begin{equation} \alpha^* = \alpha(1-\eta) + \frac{\rho a\eta }{\rho a + 1-\rho} \pm\left(2 \eta^2 + \frac{1-\rho}{(\rho a + 1-\rho)t} + \eta^{3/2}\right)\end{equation}
Hence, recalling that
\begin{equation}g_p(\alpha) \equiv -\alpha + \frac{\rho a}{\rho a + 1-\rho} = -\alpha + \frac{\alpha^pe^{c_p'(\alpha)}}{\alpha^pe^{c_p'(\alpha)} + (1-\alpha)^p}\commaeq\end{equation}
and noticing that there exists a constant $Q_\delta$ such that (for
the above range of $\alpha,\eta$, etc)
$$\frac{1-\rho}{(\rho a + 1-\rho)t}\leq \frac{Q_\delta}{t}\commaeq$$
we deduce the following bound:
\begin{multline}\Prp{[t,\alpha]}{[t,\alpha]\mapsto [(1+\eta)t,\alpha^*]\suchthat\frac{\alpha^*-\alpha}{\eta} = g_p(\alpha)\pm (2\eta + \sqrt{\eta} + Q_\eta/\eta t) \mid\elead} \\
\geq 1-e^{-V_\delta\, \eta\, t}\periodeq\end{multline} This holds
for all $\delta<\alpha<1/2-\delta$, $0<\eta<\eta_\delta$ and $t\geq
T_\delta$. Notice also that we had assumed that $\alpha t\in\N$, but
this restriction is unnecessary if we make $Q_\delta$ larger. To
summarize our conclusions, we state them in slightly modified form
as the following Lemma.

\begin{lemma}\label{lem:ODE_step} For each $\delta\in(\delta,1/2-\delta)$, there exist constants $\eta_\delta\in\R^+$, $V_\delta, Q_\delta\in\R^+$ and $T_\delta\in\N$ such that for all integers $t^*>t\geq T_\delta$ such that $0<\eta=\frac{t^*}{t}-1<\eta_\delta$,
\begin{multline}\Prp{[t,\alpha]}{[t,\alpha]\mapsto [(1+\eta)t,\alpha^*]\suchthat\frac{\alpha^*-\alpha}{\eta} = g_p(\alpha)\pm Q_\eta\left(\sqrt{\eta} + \frac{1}{t^*-t}\right) \mid\elead} \\ \geq 1-e^{-V_\delta\,(t^*-t)}\periodeq\end{multline}\end{lemma}

Notice that the constants appearing in the Lemma might be slightly
different than those appearing before it, but this is nothing but a
slight abuse of notation.

\begin{remark}\label{rem:ODE_increasing}Let us now show why $g_p(\alpha)\geq 0$ always, as stated in the introduction to this chapter. Choose $\delta<\alpha^*<\alpha<1/2-\delta$ be fixed, but also close enough to $\alpha$. Then $c_p(\alpha)\geq c_p(\alpha^*)$ always (by \thmref{escape_main}), and one can easily deduce from the reasoning proving \eqnref{ODE_goodbound} and \lemref{ODE_step} that for all $\eta>0$ fixed (but small enough)
$$\Prp{[t,\alpha]}{[t,\alpha]\mapsto [(1+\eta)t,\alpha^*]\mid \elead}\leq e^{-V_\delta\eta t}\periodeq$$
There are $\eta t$ choices for the number of balls in the first bin
at time $\eta t$, and one can easily deduce via an union bound that
$$\Prp{[t,\alpha]}{[t,\alpha]\mapsto [(1+\eta)t,\alpha^*]\mbox{ for some }\alpha^*\leq \alpha\mid \elead}\leq e^{-\ohmega{t}}\commaeq$$
(Notice that, as shown above, only
$\alpha^*>\alpha/(1+\eta)>\alpha-\eta$ need to be considered, so one
can pick $\eta$ so that $\alpha^*>\delta$ for all relevant
$\alpha^*$.) On the other hand, we know from \lemref{ODE_step} that
with very high probability $[t,\alpha]$ evolves into a state
$[(1+\eta)t,\alpha^*]$ with $\alpha^*-\alpha = g(\alpha) \eta +
\bigoh{\eta^{3/2}+ 1/ t}$. If $g(\alpha)<0$, we could pick a small
enough $\eta$ and a large enough $t$ such that $\alpha^*<\alpha$
with overwhelming probability, but this was shown to be impossible
above.\end{remark}

\subsection{Proof of \thmref{ODE_main}}\label{sub:ODE_proof}

We now have what we need to prove \thmref{ODE_main}.

\begin{proof} [of \thmref{ODE_main}] The idea of the proof is to iterate uses of \lemref{ODE_step}, which shows that, typically speaking, for any small $\eta$, the fraction $\alpha^*$ of balls after $\eta t$ time steps in bin $1$ stays close to the straight line passing through $\alpha$ with slope $g_p(\alpha)$. Since the solution $A(\cdot)=A_{p,\alpha}(\cdot)$ of the ODE also stays close to those straight lines (at least locally), this technique will give us the desired result.
Let $\alpha\in(0,1/2)$ and $0<K<T_{p,\alpha}$ be as in the statement
of the theorem. Choose some $L\in(K,T_{p,\alpha})$ and let $\delta$
be such that
\begin{equation}\delta\leq \min\left\{\alpha - \eps, \frac{1}{2} - A(L)-\eps\right\}\end{equation}
for some $0<\eps\min\{\alpha,1/2-A(L)\}$ (as discussed in the
introduction, $A(s)<1/2$ for $s<T_{p,\alpha}$, so such an $\eps$
exists).

Now recall the notation from \lemref{ODE_step}, and assume (as we
might) that $t$ satisfies
\begin{equation}\label{eq:ODE_goodt}t\geq T_\delta\commaeq \eta=\eta_t\equiv \frac{\lceil t^{1/3}\rceil}{t}<\min\{\eta_\delta,L-K\}\end{equation}
Clearly, all that \eqnref{ODE_goodt} requires is that $t$ is large
enough. Assuming it holds, there exists an integer $N_t\in\N$ such
that
\begin{equation}K\leq \eta_t N_t<L\end{equation}
and we will assume, without loss of generality, that in fact
$\eta_t\,N_t = K<L$. We also define, for convenience
\begin{equation}G^{(r)}_\delta\equiv \sup_{\delta\leq z\leq 1/2-\delta}\frac{1}{r!}\left|\frac{d^rg_p}{dz^r}(z)\right|\;\;(r\in\N)\commaeq\end{equation}
which is a finite quantity since $g$ is infinitely differentiable on
$(0,1/2)$.

Recall that we are starting the balls-in-bins process from state
$[t,\alpha]$. We will first look at the differences
\begin{equation}\label{eq:ODE_defdelta}\Delta_j\equiv |\hat\alpha(j\eta) - A(j\eta)|\commaeq\;j\in[N_t]\cup\{0\}\commaeq\end{equation}
and show that these differences remain small with high probability.
At $j=0$,
$$\hat\alpha(0) = \frac{\lceil\alpha t\rceil}{t} = A(0)\pm \frac{1}{t}$$
so the differences are small at the start. Now assume that, for some
$j\in[N_t]\cup\{0\}$, after conditioning on $\elead$,
\begin{equation}\label{eq:ODE_induction}\mbox{ With probability }\geq 1 - P_j\commaeq \\ \forall i\in[j]\cup\{0\}\commaeq\;\Delta_i\leq \gamma_i\leq \eps\periodeq\end{equation}
where $$0\leq \frac{1}{t} = \gamma_0\leq \gamma_1\leq \dots\leq
\gamma_j<\eps$$ and $P_j\in\R^+$. We will show that (again
conditioning on $\elead$)
\begin{eqnarray}\label{eq:ODE_induction2}\mbox{ With probability }\geq 1 - P_j - e^{-V_\delta \frac{t^{1/3}}{1+K}}\commaeq \\
\nonumber \forall i\in[j+1]\cup\{0\}\commaeq\;\Delta_i\leq
\gamma_i\commaeq  \mbox{where } \\ \nonumber\gamma_{j+1} \equiv
\gamma_j +(G^{(1)}_\delta+G^{(2)}_\delta)\eta^2 + Q_\delta
\left(\eta^{3/2} + \frac{1+K}{\eta t}\right)\periodeq\end{eqnarray}
(Here, $W_\delta>0$ is a constant depending only on $\delta$). To
prove this, let us condition on a value of $\hat\alpha_t(\eta j)$
that is compatible with the event described in
\eqnref{ODE_induction}. This means that
\begin{equation}\label{eq:ODE_gooda}\hat\alpha_t(\eta j)= \alpha_j \mbox{ with } |A(\eta j)-\alpha_j| \leq \gamma_j\periodeq\end{equation}
In this case, since $\gamma_j<\eps$ and $A$ is increasing,
\begin{equation} \delta< \alpha - \eps = A(0)-\eps \leq \alpha_j\leq A(L) + \eps <\frac{1}{2}-\delta\periodeq\end{equation}
Using the Markov Property of the balls-in-bins process shows that
\begin{equation}\label{eq:ODE_toacabando}\Prp{[t,\alpha]}{\Delta_{j+1}\leq \gamma_{j+1}\mid\elead\commaeq \hat\alpha_t(\eta j)=\alpha_j} = \Prp{[(1+j\eta) t,\alpha_j]}{\Delta_{j+1}\leq \gamma_{j+1}\mid\elead}\end{equation}
since $\hat\alpha_t(s)$ is the number of balls in bin $1$ at time
$\lceil st\rceil$ (i.e. when there are $t+\lceil st\rceil$ balls in
the system) and in the present case $\eta t\in\N$. To evaluate the
latter probability, notice first that
\begin{eqnarray} & & \left|A((j+1)\eta) - A(j\eta)  - \eta\,A'(j\eta) \right|\\ & = &  \left|A((j+1)\eta) - A(j\eta)  - \eta g_p(A(j\eta)) \right|\\
&\leq & G^{(2)}_\delta \eta ^2\periodeq\end{eqnarray} Moreover, by
\eqnref{ODE_gooda} and the choice of $t\geq T_\delta$ one can apply
\lemref{ODE_step} with $(1+j\eta) t$ replacing $t$ and
$\eta/(1+j\eta)$ replacing $\eta$ to deduce that, conditioned on
$\hat\alpha_t(j\eta)=\alpha_j$ as above, the probability that
\begin{equation} \left|\hat\alpha_t((j+1)\eta) - \hat\alpha_t(j\eta) - \eta g_p(\hat\alpha_t(j\eta)) \right|\leq Q_\delta \left(\eta^{3/2} + \frac{1+j\eta}{\eta t}\right)\end{equation}
is at least $1-e^{-V_\delta \eta t/(1+j\eta)}$. When the two
previous equations hold,
\begin{multline}\label{eq:ODE_evolvesright}|\hat\alpha_t(s_0+\eta) - A(s_0+\eta) | \leq |\hat\alpha_t(s) - A(s_0)| + \eta |g_p(\hat\alpha_t(s_0)) - g_p(A(s_0))| \\ + G^{(2)}_\delta\eta^2  + Q_\delta \left(\eta^{3/2} + \frac{1}{\eta t}\right) \leq \gamma_j + (G^{(1)}_\delta+G^{(2)}_\delta)\eta^2 + Q_\delta \left(\eta^{3/2} + \frac{1+K}{\eta t}\right)\periodeq\end{multline}
Thus, for any $\alpha_j$ compatible with $\Delta_j\leq \gamma_j$,
one has that
\begin{equation}\label{eq:ODE_toacabando2}\Prp{[t,\alpha]}{\Delta_{j+1}\leq \gamma_{j+1}\mid\elead\commaeq \hat\alpha_t(\eta j)=\alpha_j}\geq 1 - e^{-V_\delta \eta t/(1+K)}\commaeq\end{equation}
from which \eqnref{ODE_induction2} immediately follows.

Now notice that if
$$\gamma_N \equiv N \left\{(G^{(1)}_\delta+G^{(2)}_\delta)\eta^2 + Q_\delta \left(\eta^{3/2} + \frac{1+K}{\eta t} + \frac{1}{t}\right)\right\}\leq \eps$$
then one can use \eqnref{ODE_induction} and \eqnref{ODE_induction2}
repeatedly to deduce that
\begin{equation}\label{eq:ODE_discretepoints}\Prp{[t,\alpha]}{\forall j\in[N]\cup\{0\}\commaeq
\Delta_j\leq \gamma_N\mid \elead}\geq 1 - N\,e^{-V_\delta \eta
t/(1+K)}\periodeq\end{equation} But notice that $N = K/\eta$, so a
simple calculation shows that
$$\mbox{for all large enough $t$}\commaeq \gamma_N \leq W_\delta\, t^{-1/3}$$
with $W_\delta\in\R^+$ depending only on $\delta$. Hence
$\gamma_N\leq \eps$ for all large enough $t$, and for such $t$
\eqnref{ODE_discretepoints} holds. Finally, one can easily show that
in the event described by \eqnref{ODE_discretepoints},
\begin{equation}\forall j\in[N-1]\cup\{0\}\commaeq\forall s\in(0,\eta)\commaeq |\hat\alpha_t(j\eta+s) - A(s)|\leq 2\eta + \gamma_N =\bigoh{t^{-2/3}}\commaeq t\gg 1\periodeq\end{equation}
Hence, \eqnref{ODE_discretepoints} actually implies that for all
large enough $t$,
\begin{equation}\label{eq:ODE_allpoints}\Prp{[t,\alpha]}{\sup\limits_{s\in[0,K]} |\hat\alpha_t(s) - A(s)|\leq W_\delta t^{-1/3}\mid \elead}\geq 1 - \bigoh{t^{2/3}}\,e^{-V_\delta
t^{1/3}/(1+K)}\commaeq\end{equation} for a possibly larger
$W_\delta$. Since $\delta$ is ultimately defined in terms of
$\alpha$ and $K$, \eqnref{ODE_allpoints} implies the
Theorem.\end{proof}

\subsection{Proof of \lemref{ODE_Chernoff}}\label{sub:ODE_Chernoff}
To conclude the chapter, we prove \lemref{ODE_Chernoff}.
\begin{proof} [of \lemref{ODE_Chernoff}] Notice first that, if $0<n<m$
\begin{equation}\label{eq:ratiob}\frac{b(n)}{b(n+1)} = \frac{n+1}{m-n} \frac{1-p}{ap}\end{equation}
Notice that $x\mapsto (x+1/m)/(1-x) = (x+1/m)\sum_{\ell\geq
1}x^\ell$ is an increasing function of $x$ that is equal to $1/m<1$
at $x=0$ and goes to $+\infty$ as $x\to 1$. Hence, if
$$x_0= \frac{\frac{pa}{1-p}-\frac{1}{m}}{1 +\frac{pa}{1-p}} = \frac{pa-\frac{1-p}{m}}{pa + (1-p)}$$
then
\begin{equation}\label{eq:defx}\forall x\in[0,1)\;\;  \left\{\begin{array}{lllllll}\frac{x+\frac{1}{m}}{1-x} \frac{1-p}{ap} &>& 1&\iff & x&>&x_0 \\
\frac{x+\frac{1}{m}}{1-x} \frac{1-p}{ap} &=& 1&\iff& x&=&x_0 \\
\frac{x+\frac{1}{m}}{1-x} \frac{1-p}{ap} &<& 1&\iff &
x&<&x_0\end{array}\right.\end{equation} As a result, if we let
$n_0\equiv \lceil x_0m \rceil$ (which is the same definition as in
the statement of the lemma), we have that (using \eqnref{ratiob})

\begin{equation*}\forall 0<j\leq m-n_0 \;\;\frac{b(n_0)}{b(n_0+j)} = \prod_{i=1}^{j}\frac{b(n_0+i-1)}{b(n_0+i)} > 1\end{equation*}

and similarly
\begin{equation*}\forall 0<j\leq n_0 \;\;\frac{b(n_0)}{b(n_0-j)} = \prod_{i=1}^{j}\frac{b(n_0-i+1)}{b(n_0-i)} < 1\end{equation*}

This proves the first two items in the lemma. As for the last one,
it suffices to show that for all $j>K\sqrt{m}$
$$b(n_0+j),b(n_0-j)\leq b(n_0)e^{-\frac{K}{2(1-x^*)}}$$
We only prove the first inequality; the proof of the second is
almost identical. As before, we write (using \eqnref{ratiob})

$$b(n_0+j) = b(n_0)\prod_{i=1}^{j}\frac{b(n_0+i)}{b(n_0+i-1)} = b(n_0) \left(\prod_{i=1}^{j}\frac{m-(n_0+i)}{n_0+i+1}\right) \left(\frac{pa}{(1-p)}\right)^j$$
Now notice that (using the definition of $x_0$ and $n_0$)

\begin{eqnarray*}\left(\prod_{i=1}^{j}\frac{m-(n_0+i)}{n_0+i+1}\right) &<& \left(\frac{m-n_0}{n_0+1}\right)^j \prod_{i=1}^{j}\left(1-\frac{i}{m-n_0}\right) \\
&\leq& \left(\frac{m-n_0}{n_0+1}\right)^j \exp\left(-\sum_{i=1}^j\frac{i}{m-n_0}\right) \\
&= &  \left(\frac{1-x}{x+1/m}\right)^j \exp(-\frac{j(j+1)}{2(1-x)m})
\end{eqnarray*}

where $x\equiv n_0/m$ is bigger than $x_0$. As a result of
\eqnref{defx}

$$ \frac{1-x}{x+1/m} \frac{ap}{1-p}<1$$

and putting this together with the previous inequalities

$$\frac{b(n_0+j)}{ b(n_0) } <\left(\frac{pa}{(1-p)}\right)^j\left(\frac{1-x}{x+1/m}\right)^j \exp(-\frac{j(j+1)}{2(1-x)m}) \leq \exp(-\frac{j(j+1)}{2(1-x)m})$$

which, together with the fact that $j>K\sqrt{m}$, finishes the
proof.\end{proof}

\section{Open problems}\label{sec:final}

The results proven here only apply to feedback functions $f(x)=x^p$.
However, we have been able to prove other results for more general
functions; see \cite{MitzenmacherOS04,EuBrown,EuOnset} for several
examples. It would be interesting to see extensions of the present
work to those other feedback functions as well.

Another open problem is to determine the asymptotic behaviour of the
ordinary differential equation in \thmref{ODE_main}, especially
whether the solution blows up in finite time (i.e.
$T_{p,\alpha}<+\infty$). We conjecture that this is {\em not} the
case, but a proof would require a careful analysis of the ODE.

\bibliography{escape}
\bibliographystyle{plain}

\end{document}